\documentclass[12pt]{amsart}

\usepackage{amssymb,latexsym,amscd}

\usepackage{amsthm,amsmath}

\usepackage{xypic}

\newcommand{\cK}{\mathcal{K}}

\newcommand{\Sym}{\mathrm{Sym}}

\newcommand{\lra}{\longrightarrow}

\newcommand{\ra}{\rightarrow}

\newcommand{\M}{\mathcal{M}}

\newcommand{\PP}{\mathbb{P}}

\newcommand{\ZZ}{\mathbb{Z}}

\newcommand{\CC}{\mathbb{C}}

\newcommand{\so}{\mathrm{SO}}

\theoremstyle{plain}

\newtheorem{thm}{Theorem}[section]

\newtheorem{lem}[thm]{Lemma}

\newtheorem{prop}[thm]{Proposition}

\newtheorem{cor}[thm]{Corollary}

\theoremstyle{remark}

\newtheorem{rem}[thm]{Remark}

\begin{document}

\title[]{On Theta Functions of Order Four}

\author{Yaacov Kopeliovich}

\author{Christian Pauly} 

\author{Olivier Serman}

\address{MEAG NY \\ 540 Madison Avenue 6th floor \\ New York \\ NY 10022 \\ USA}
\email{ykopeliovich@yahoo.com}

\address{D\'epartement de Math\'ematiques \\ Universit\'e de Montpellier II - Case Courrier 051 \\ Place Eug\`ene Bataillon \\ 34095 Montpellier Cedex 5 \\ France}
\email{pauly@math.univ-montp2.fr}

\address{Laboratoire Paul Painlev\'e \\ UMR 8524 du CNRS \\ Universit\'e des Sciences et Technologies de Lille \\ Cit\'e scientifique \\ 59655 Villeneuve d'Ascq Cedex \\ France}

\email{Olivier.Serman@math.univ-lille1.fr}

%\thanks{}

%\keywords{}


\begin{abstract}
We prove that the fourth powers of theta functions with even characteristics form a basis of the space $H^0(A,\mathcal O_A(4\vartheta))_+$ of even theta functions of order four on a principally polarized Abelian variety $(A,\vartheta)$ 
without vanishing theta-null.
\end{abstract}

\maketitle

%%%%%%%%%%%%%%%%%%%%%%%%%%%%%%%%%%%%%%%%%%%%%%%%%%%%%%%%%%

\section{Introduction}

Let $(A,\vartheta)$ be a principally polarized Abelian variety of dimension $g$ defined over an algebraically closed field
$k$ of characteristic different from two.
By Mumford's algebraic theory of theta functions 
\cite{M} (see also the appendix of \cite{B}) there is a canonical bijection between the 
set of symmetric effective divisors 
$\Theta_\kappa \subset A$ representing the principal polarization and the set $\cK$ of theta-characteristics $\kappa$ of $A$.  We denote by $L$ the unique symmetric line bundle 
over $A$ representing $2\vartheta$ and such that the linear system $|L|$ contains the
divisors $2\Theta_\kappa$ for all $\kappa \in \cK$.  
The space $H^0(A,L^2)$ of theta functions of order four decomposes
under the natural involution of $A$ into $\pm$-eigenspaces $H^0(A,L^2)_+$ and $H^0(A,L^2)_-$ of dimensions 
$$ d_+ = 2^{g-1}(2^g+1) \qquad \text{and}  \qquad d_- = 2^{g-1}(2^g-1). $$
We note that $d_+$ and $d_-$ are equal to the number of even and odd theta-characteristics of $A$.
Moreover any even theta-characteristic 
$\kappa_0 \in \cK_+$ decomposes according to the values $\pm 1$ taken by 
$\kappa_0$ the set of $2$-torsion points
$A[2]$ into a union of two subsets 
$A[2]_+$ and $A[2]_-$ of cardinality $d_+$ and $d_-$. Note that $0 \in A[2]_+$. We denote by 
$$ \varphi : A \longrightarrow |L^2|_+^* = \PP^{d_+ -1} $$
the morphism induced by the linear system $|L^2|_+$. Then our main result is the following

\begin{thm}\label{projectivebasis}
For any even theta-characteristic $\kappa_0 \in \cK_+$,
the $d_+$ points 
$$ \varphi(a) \in |L^2|^*_+ = \PP^{d_+ -1} \qquad \text{with} \qquad a \in A[2]_+$$
form a projective basis if and only if $A$ has no vanishing theta-null.
\end{thm}

There is also a ``dual'' version of this result. 

\begin{thm}\label{functions}
The $d_+$ divisors 
$$ 4 \Theta_\kappa \in  |L^2|_+ \qquad \text{with} \qquad \kappa \in \cK_+$$
form a projective basis (or, in other words, the fourth powers of theta functions
with even characteristics form a basis of $H^0(A,L^2)_+$) if and only if $A$ has no vanishing theta-null.
\end{thm}

Although of independent interest, these results also have some consequences on {\em
generalized} theta functions (see \cite{S} for details): let $\M_{\so_r}^+$ denote the
moduli space of topologically trivial oriented orthogonal bundles of rank $r$ over a smooth 
projective complex curve $C$ of genus $g \geq 2$ and let $\mathcal L_{\so_r}^+$ denote the determinant
line bundle over $\M_{\so_r}^+$. It was shown in \cite{Borth} that the linear system 
$|\mathcal L_{\so_r}^+|$ identifies canonically with the dual $|r \Theta|^*_+$ of even
theta functions of order $r$ over the Picard variety parametrizing degree $g-1$ line
bundles over $C$. In particular, for $r=4$ we observe that 
$\dim H^0(\M_{\so_4}^+, \mathcal L_{\so_4}^+) = d_+$. Moreover we can associate to any
even theta-characteristic $\kappa \in \cK_+$ a proper effective divisor $\Delta_\kappa \in 
|\mathcal L_{\so_4}^+|$ with support
$$ \{ E  \in  \M_{\so_4}^+ \ | \ \dim  H^0(E\otimes \kappa) \neq 0\}.$$
By pulling back the $d_+$ divisors $\Delta_\kappa$ to the Jacobian $\mathrm{Jac}(C)$ under
the map $\mathrm{Jac}(C) \ra \M_{\so_4}^+$, $ L \mapsto L \oplus L^{-1} \oplus L \oplus L^{-1}$,
we deduce from Theorem 1.2 the following

\begin{cor}
Assume that $C$ has no vanishing theta-null. Then the $d_+$ divisors $\Delta_\kappa \in
|\mathcal L_{\so_4}^+|$  for $\kappa \in \cK_+$ form a projective basis.
\end{cor}

Finally, we mention that for general $r$ the divisors $\Delta_\kappa \in
|\mathcal L_{\so_r}^+|$ --- more precisely, some relatives of $\Delta_\kappa$ on the
moduli stack of $\mathrm{Spin}_r$-bundles over $C$ --- have recently been studied in connection with the 
strange duality of generalized theta functions \cite{Be}.

\bigskip

We would like to thank A. Beauville and S. Grushevsky for their interest in these questions. The 
second author was partially supported by the Ministerio de Educaci\'on y Ciencia (Spain) through the grant
SAB2006-0022.

\section{Review of algebraic theory of theta functions}

In this section we recall the results on theta functions of order four, which we need in the
proof of the two main theorems. We refer to \cite{M} and \cite{B} for details and proofs.

\medskip

We recall that a theta-characteristic is a map $\kappa : A[2] \rightarrow \{ \pm 1 \}$ such that
$\kappa(a+b) = \kappa(a) \kappa(b) \langle a,b \rangle$, where $\langle \ , \ \rangle$ is the
symplectic Weil form on $A[2]$. The group $A[2]$ acts transitively on the set $\cK$ of
theta-characteristics by the formula 
\begin{equation} \label{thetacar}
(a \cdot \kappa)(b) = \langle a,b \rangle \kappa(b).
\end{equation}
We say that $\kappa$ is even if $\kappa$ takes $d_\pm$ times the value $\pm 1$. The
set of even theta-characteristics is denoted by $\cK_+$.

\medskip

The Heisenberg group $H$ associated to the line bundle $L$ consists of all pairs 
$\alpha = (a, \varphi)$ with $a \in A[2]$ and $\varphi : T^*_a L 
\stackrel{\sim}{\rightarrow} L$ an isomorphism (where $T_a \colon A \to A$ is the translation by $a$), and fits into an exact sequence
$$ 0 \lra k^* \lra H \stackrel{p}{\longrightarrow} A[2] \lra 0, \qquad
\text{with} \ p(\alpha) = a.$$
For any $\kappa \in \cK$ we consider  the character $\chi_\kappa : H
\ra k^*$ of weight $2$ defined by $\chi_\kappa (\alpha) = || \alpha || \kappa(a)$, with
$|| \alpha || = \alpha^2 \in k^*.$ The main result we need in this note is

\begin{prop}[\cite{B} Proposition A.8]  \label{recall}
The vector space $H^0(A, L^2)_+$ decomposes as an $H$-module into a direct sum
$\bigoplus_{\kappa \in \cK_+} W_\kappa$ with $\dim W_\kappa$ = 1 and $W_\kappa$
is the character space associated to the character $\chi_\kappa$. The zero divisor of any
nonzero section  $s_\kappa \in W_\kappa$ equals $2_A^* \Theta_\kappa$, where $2_A$ is the
duplication map of $A$.
\end{prop}

\section{The matrix $M$}

We say that a pair $(a_1, a_2) \in (\ZZ/ 2 \ZZ)^g \times (\ZZ/ 2 \ZZ)^g$ is even if $a_1 \cdot a_2 = 0$. Let $M$ denote the
$d_+ \times d_+$ matrix with lines and columns indexed by even pairs $(a_1, a_2)$ and $(b_1, b_2)$ and
with entries
$$ (-1)^{a_1 \cdot b_2 + a_2 \cdot b_1}.$$

\noindent We will use the fact that $M$ is invertible, which has already been proved in \cite{fay} Lemma 1.1. For the convenience of the reader, and for later use, we recall here the computations.

\begin{lem}\label{keyidentity}
For fixed $(a_1, a_2) \in (\ZZ/ 2 \ZZ)^{g} \times (\ZZ/ 2 \ZZ)^g$, we have the equality
$$ \sum_{b_1, b_2 \in (\ZZ/ 2 \ZZ)^{g}, \atop{b_1 \cdot b_2 = 0}}
(-1)^{a_1 \cdot b_2 + a_2 \cdot b_1} =  
\begin{cases}
d_+ & \text{if } \ a_1 = a_2 = 0, \\
(-1)^{a_1 \cdot a_2} 2^{g-1} & \text{else.}
\end{cases}$$
\end{lem}

\begin{proof}
If $a_1=a_2=0$ there is nothing to prove. Otherwise we may assume that $a_1 \neq 0$, and the result follows from the easy equality
$$\sum_{b_1 \in (\ZZ/2\ZZ)^g, \atop{b_1 \cdot b_2=0}} (-1)^{a_1 \cdot b_1}=\begin{cases} 2^{g-1} & \text {if } b_2 = a_1 \\ 0 & \text{else. } \end{cases}$$
\end{proof}

\begin{prop}\label{inverse}
The matrix $M$ is invertible.
\end{prop}

\begin{proof} \label{Minv}

Using Lemma \ref{keyidentity}, we easily check that the inverse of $M$ is given by
$$M^{-1}=\genfrac{}{}{}{}{1}{2^{2g-1}}(M-2^{g-1}I).$$

\end{proof}

\section{Proof of Theorem \ref{projectivebasis}}

Consider the linear map given by evaluating the $d_+$ theta functions 
$s_\kappa \in H^0(A,L^2)_+$ at the $d_+$ points $ a \in A[2]_+$
$$ \mathrm{ev} : H^0(A, L^2)_+ \longrightarrow \bigoplus_{a \in A[2]_+} L^2_a.$$
\noindent It is clear that the assertion of Theorem \ref{projectivebasis} is equivalent to $\mathrm{ev}$
being an isomorphism. First note that saying that $A$ has a vanishing theta-null means that there exists a
section $s_\kappa$, with $\kappa \in \cK_+$, which vanishes at every point $a\in A[2]$, i.e. $s_\kappa \in
\ker \mathrm{ev}$. On the other hand we will show that, if $A$ has no vanishing theta-null, $\mathrm{ev}$ is given after suitable normalization by the matrix $M$. The theorem then follows from Proposition \ref{inverse}.

\medskip

We consider for each $\kappa \in \cK_+$ a nonzero section $s_\kappa \in H^0(A,L^2)_+$ in the
one-dimensional $\chi_\kappa$-character space $W_\kappa$, and for each $a \in A[2]$ an
isomorphism $\phi_a : L^2_a \stackrel{\sim}{\rightarrow} k$. Since we assume that $A$
has no vanishing theta-null, we have $\phi_a(s_\kappa(a)) \not= 0$ for all 
$\kappa \in \cK_+$ and all $a \in A[2]$. The quotient
$$ \mu(a, \kappa, \kappa') = 
\frac{\phi_0(s_\kappa(0)) \cdot \phi_a(s_{\kappa'}(a))}{\phi_a(s_\kappa(a)) \cdot \phi_0(s_{\kappa'}(0))}, 
\qquad a \in A[2], \ \kappa, \kappa' \in \cK_+$$
does not depend on the choice of the sections $s_\kappa, s_{\kappa'}$ and the isomorphisms
$\phi_a$. Given $a \in A[2]$ we choose an $\alpha = (a, \varphi) \in H$ such that 
$\varphi: T_a^* L \stackrel{\sim}{\rightarrow} L$ preserves the isomorphisms
$\phi_0$ and $\phi_a$, which is equivalent to the equality $\phi_0( (\alpha. s_\kappa)(0)) =
\phi_a (s_\kappa(a))$ for all $\kappa \in \cK_+$. On the other hand $\alpha. s_\kappa 
= \chi_\kappa (\alpha) s_\kappa =  ||\alpha || \kappa(a) s_\kappa$, hence
$\frac{\phi_a(s_\kappa(a))}{\phi_0(s_\kappa(0))} = || \alpha || \kappa(a)$. 
Therefore we obtain the formula 
\begin{equation} \label{eqcoeff}
\mu(a, \kappa, \kappa') =  \kappa(a) \kappa'(a).
\end{equation}
In order to obtain the matrix $M$ we normalize as follows: consider a section $s_{\kappa_0}$ corresponding to the fixed even theta-characteristic $\kappa_0$, and choose the isomorphisms $\phi_a$ such that $\phi_a( s_{\kappa_0}(a)) = 1$ for all $a \in A[2]_+$. Then we choose the sections $s_\kappa$ such that $\phi_0(s_\kappa(0)) = 1$. Any $\kappa \in \cK$ can be written
$\kappa = b \cdot \kappa_0$ for a unique  $b \in A[2]$, and $\kappa \in \cK_+$ if and only if $b \in A[2]_+$ by \cite{B} formula (1) page 279. 
Using formulae \eqref{eqcoeff} and  \eqref{thetacar} we obtain the equalities
$$\phi_a(s_{b \cdot \kappa_0}(a)) = \mu ( a, \kappa_0, b \cdot \kappa_0 ) = 
\kappa_0(a) (b \cdot \kappa_0)(a) = \langle a,b \rangle. $$ 
We now choose a level-$2$ structure  $\lambda : A[2] \stackrel{\sim}{\rightarrow}
(\ZZ/ 2 \ZZ)^g \times (\ZZ/ 2 \ZZ)^g$ which maps the set $A[2]_+$ bijectively onto the
set of even pairs --- this is possible, since $\mathrm{Sp}(2g, \ZZ/2\ZZ)$ acts transitively on the
set of even theta-characteristics. Then we can write $\langle a,b \rangle = 
(-1)^{a_1 \cdot b_2 + a_2 \cdot b_1}$ with $\lambda(a) = (a_1,a_2)$ and
$\lambda(b) = (b_1,b_2)$. This finishes the proof.

\section{Proof of Theorem \ref{functions}}
\label{sectionfunctions}
Let $\kappa_0$ be a theta-characteristic of $A$. Recall from \cite{B} that the morphism $\delta _{\kappa_0} \colon A \lra |L|$ which sends $a \in A$ to the divisor $T_a^\ast \Theta_{\kappa_0} + T_{-a}^\ast \Theta_{\kappa_0}$ fits into the commutative diagram

$$\xymatrix@R=20pt@C=50pt{
 & |L|^{\ast} \ar^\wr[dd]\\
A \ar_{\delta_{\kappa_0}}[dr] \ar^{\varphi_L}[ur] & \\
& |L|,
}$$
\noindent where $\varphi_L$ is the morphism defined by the complete linear system $|L|$. The isomorphism 
between $|L|^\ast$ and $|L|$ is given by any nonzero element in the one-dimensional character space in $H^0(A,L) \otimes H^0(A,L)$ associated to $\chi_{\kappa_0}$. Now, if $A$ has no vanishing theta-null, we know from \cite{B} Proposition A.9 that the multiplication map
$$\Sym^2 H^0(A,L) \lra H^0(A,L^2)_+$$
\noindent is bijective. It follows that the duplication morphism $|L| \lra |L^2|_+, D \longmapsto 2D$ is identified, through the isomorphism $|L|^\ast \buildrel\sim\over\lra |L|$, with the $2$-uple embedding of $|L|^\ast$, and that the morphism $\varphi \colon A \lra |L^2|_+^\ast$ is the composite of $\varphi_L$ with this $2$-uple embedding.

But we know from Theorem \ref{projectivebasis} that the images $\varphi(a)$ of the $d_+$ points $a \in A[2]_+$ form a projective basis of $|L^2|_+^\ast$. This implies that the points $2\delta_{\kappa_0}(a)=4\Theta_{a \cdot \kappa_0}$ form a projective basis of $|L^2|_+$, which finishes the proof.

\begin{rem}
 The previous considerations show that our result is also equivalent to the following one: if $A$ has no vanishing theta-null, there is no quadric hypersurface in $|L|$ containing the $d_+$ points defined by $2\Theta_\kappa$, with $\kappa \in \cK_+$.
\end{rem}

\begin{rem} 
It easily follows from the previous proof that the codimension of the linear span of the fourth powers of even
theta functions in $H^0(A, L^2)_+$ equals the number of vanishing theta-nulls of $A$.
\end{rem}

\section{An analytic proof of Theorem 1.2}

When $k=\CC$, we can give a short analytic proof of Theorem 1.2 by using Riemann's quartic addition theorem to express the fourth powers of the theta functions $\theta_\kappa$ in terms of the functions $2_A^\ast\theta_\kappa$.

In this case, the Abelian variety $A$ is a quotient $\CC^g / \Gamma_\tau$ of a $g$-dimensional vector space by a lattice $\Gamma_\tau=\ZZ^g \oplus \tau \ZZ^g$ for some $\tau$ in the Siegel upper half-plane.
Let us choose a level-$2$ structure to identify the set of theta-characteristics with $(\ZZ/2\ZZ)^g \times (\ZZ/2\ZZ)^g$, in such a way that even theta-characteristics $\kappa$ correspond to even pairs $(a_1,a_2)$ (see \cite{B} A.6). We want to prove that the space $H^0(A,L^2)_+$ of even theta functions of order four is spanned by the fourth powers of the theta functions with even characteristics $a_1,a_2 \in (\ZZ/2\ZZ)^g$
$$\displaystyle{\theta{\genfrac{[}{]}{0pt}{}{a_1}{a_2}}(z) =\sum_{m \in \ZZ^g} \exp \pi i \bigl( {}^t(m+\genfrac{}{}{}{}{1}{2} \tilde{a_1}) \tau (m+\genfrac{}{}{}{}{1}{2} \tilde{a_1}) + 2 \, {}^t(m+\genfrac{}{}{}{}{1}{2}\tilde{a_1})(z+\genfrac{}{}{}{}{1}{2}\tilde{a_2})\bigr)},$$
\noindent where $\tilde{a_1}$ and $\tilde{a_2}$ denote some representatives in $\ZZ^g$ of $a_1$ and $a_2$ 
respectively.

Riemann's quartic addition theorem (see \cite{tata} II.6 formula $(\text R _{\text{ch}})$) implies that, for any pair $(a_1,a_2)$,
$$\bigl(\theta {\genfrac{[}{]}{0pt}{}{a_1}{a_2}}(z)\bigr)^4=2^{-g} \sum_{b_1,b_2 \in (\ZZ/2\ZZ)^g} (-1)^{a_1\cdot b_2 + a_2\cdot b_1} \bigl(\theta{\genfrac{[}{]}{0pt}{}{b_1}{b_2}}(0)\bigr)^3\theta{\genfrac{[}{]}{0pt}{}{b_1}{b_2}}(2z).$$

We now invert this formula as follows. Let us sum over all even pairs $(a_1,a_2)$
the previous relations multiplied by the factor $(-1)^{a_1\cdot c_2 + a_2 \cdot c_1}$, where $(c_1,c_2)$ is a fixed even pair. Using Lemma \ref{keyidentity}, we obtain
\begin{multline*}2 \sum_{a_1, a_2 \in (\ZZ/2\ZZ)^g, \atop{a_1 \cdot a_2=0}} (-1)^{a_1 \cdot c_2 + a_2 \cdot c_1} \bigl(\theta\genfrac{[}{]}{0pt}{}{a_1}{a_2}(z)\bigr)^4 = \\ 2^g\bigl(\theta{\genfrac{[}{]}{0pt}{}{c_1}{c_2}}(0)\bigr)^3\theta{\genfrac{[}{]}{0pt}{}{c_1}{c_2}}(2z)+\sum_{b_1, b_2 \in (\ZZ/2\ZZ)^g} (-1)^{(b_2+c_2) \cdot (b_1 + c_1)}\bigl(\theta{\genfrac{[}{]}{0pt}{}{b_1}{b_2}}(0)\bigr)^3\theta{\genfrac{[}{]}{0pt}{}{b_1}{b_2}}(2z).
\end{multline*}
\noindent In the last sum, the terms corresponding to odd pairs $(b_1,b_2)$ vanish, because $\displaystyle{\theta\genfrac{[}{]}{0pt}{}{b_1}{b_2}(0)=0}$. Since $(c_1,c_2)$ is an even pair, this sum is equal to $\displaystyle{2^g \bigl(\theta\genfrac{[}{]}{0pt}{}{c_1}{c_2}(z)\bigr)^4}$, again by Riemann's quartic relation. We finally get, for every even pair $(c_1,c_2)$,
$$2^g\bigl(\theta{\genfrac{[}{]}{0pt}{}{c_1}{c_2}}(0)\bigr)^3\theta{\genfrac{[}{]}{0pt}{}{c_1}{c_2}}(2z) = - 2^g \bigl(\theta\genfrac{[}{]}{0pt}{}{c_1}{c_2}(z)\bigr)^4 + 2 \sum_{a_1, a_2 \in (\ZZ/2\ZZ)^g, \atop{a_1 \cdot a_2=0}} (-1)^{a_1 \cdot c_2 + a_2 \cdot c_1} \bigl(\theta\genfrac{[}{]}{0pt}{}{a_1}{a_2}(z)\bigr)^4.$$
\noindent But we have already recalled in Proposition \ref{recall} that the functions $\displaystyle{z \mapsto \theta{\genfrac{[}{]}{0pt}{}{c_1}{c_2}}(2z)}$ with $c_1 \cdot c_2=0$ span $H^0(A,L^2)_+$. This shows that the functions $\displaystyle{\bigl(\theta {\genfrac{[}{]}{0pt}{}{a_1}{a_2}}(z)\bigr)^4}$ with $a_1 \cdot a_2=0$ span this vector space if and only if $A$ does not have any vanishing theta-null.

\begin{rem}
It is possible to make this analytic proof work in the algebraic set-up \cite{M} of
theta functions over any algebraically closed field of characteristic different from
two. However, its algebraic version would be longer and more technical than the proof we gave in
section 5. 
\end{rem}

%It is possible to work out a similar proof in the algebraic set-up
%\cite{M} of theta functions over any algebraically closed field of characteristic
%different from 2, but this algebraic proof is much more technical and
%lengthy. For that reason we only give here the case $k= \CC$, where theta
%functions are easier to handle with, and refer to section \ref{sectionfunctions} for a short
%algebraic proof

\end{document}